\documentclass[12pt]{article}
\usepackage{amsmath}
\usepackage{amssymb}
\usepackage{graphics}
\usepackage{amsthm}
\theoremstyle{plain}

\newtheorem{thm}{Theorem}[subsection]  
\newtheorem{lem}[thm]{Lemma}           
\newtheorem{defin}[thm]{Definition}    

\newcommand{\intmean}{{-\hspace{-12.5pt}\int}}
\makeatletter
\def\theequation{\thesection.\@arabic\c@equation}
\def\thethm{\thesection.\@arabic\c@thm}
\def\thelem{\thesection.\@arabic\c@thm}
\def\thecrlr{\thesection.\@arabic\c@thm}
\def\theprp{\thesection.\@arabic\c@thm}
\def\therem{\thesection.\@arabic\c@thm}
\makeatother

\renewcommand{\epsilon}{\varepsilon}

\newcommand{\R}{{\mathbb R}}
\date{}

\begin{document}
\title{INTERIOR ESTIMATES IN CAMPANATO SPACES \\ RELATED TO
QUADRATIC FUNCTIONALS}

\author{Maria Alessandra Ragusa 
\& Atsushi Tachikawa
}

\maketitle
\begin{abstract}
In this paper we obtain interior estimates in Campanato spaces for
the derivatives of the minimizers of quadratic functionals.
\end{abstract}
\small
\begin{quote}
\end{quote}
\normalsize
\section{Introduction and Preliminary Tools}

The papers concerned with the regularity problem almost always  have as a common starting point the Euler's
equation related to a generic functional $I.$ In the paper by Giaquinta and
Giusti \cite{8}, the authors investigate the
h\"older continuity of the minima working directly with the
functional $I$ instead of Euler's equation.
 
In the present paper, following the method in \cite{8}, we have studied regularity properties of
the minima of variational integrals of the type:
$$
{\cal A}(u ,\Omega)\,=\, \int f(x,u,Du)dx,$$
$\Omega \subset \R^n, n \geq 3$ is a bounded opens set $u:\Omega \to
\R^N ,$ $u(x)=(u^1(x), \ldots, u^N (x)),$ $Du=D_\alpha u^i ,$
$D_\alpha = \frac{\partial}{\partial x^\alpha} ,$ 
$i=1,\ldots ,N, $ $\alpha=1,\ldots ,n,$  
$$
f:\Omega \times \R^N \times \R^{nN} \to \R
$$
defined by
$$
f(x,u,Du)\,=\, A^{\alpha \beta}_{i j}(x,u) D_\alpha u^i D_\beta u^j
\,+\, g(x,u,Du).
$$

Let us now give the following definitions useful in the sequel.

\begin{defin}\label{def2.1}(see \cite{K-J-F},  \cite{Ragusa1}).
Let $1\leq p<\infty, 0 \leq \lambda <n.$ 

By $L^{p,\lambda} (\Omega )$ we denote the linear space of functions
$f \in L^p(\Omega)$ such that
$$
\| f \|_{L^{p,\lambda} (\Omega)} =\biggl\{
\sup_{\underset{x \in \Omega}{0<\rho < \mathrm{diam}\, \Omega }} 
\rho^{ -\lambda } 
\int\limits_{\Omega \cap B(x,\rho)} |f (y)|^p dy\biggr\} < + \infty
$$
where $B(x,\rho)$   
ranges in the class of  the balls of $\R^n$ of radius $\rho$ around $x.$
\end{defin}
We have that $\|f \|_{L^{p,\lambda} (\Omega)} $ is a
norm respect to which $L^{p,\lambda} (\Omega)$ is a Banach space, and
also that 
$$
\|f\|_{L^{p} (\Omega)} \leq ( \mathrm{diam}
\, \Omega)^\frac{\lambda}{p} \|f \|_{L^{p,\lambda} (\Omega)}.
$$
Before the definition of the Campanato spaces let us define 
 $f_{{B(x,\rho)}}$ as the integral average
$$
f_{B(x,\rho)}= \frac{1}{|_{B(x,\rho)}|} \int\limits_{B(x,\rho)} f(y)dy
$$ 
of the function $f(x)$ over the balls ${B(x,\rho)}$ of ${\R}^n.$
When no confusion may arise, we will write $f_{B_\rho}$ or $f_\rho$ instead of $f_{B(x,\rho)}.$

\begin{defin}\label{def2.2}(see e. g. \cite{1}, \cite{4}).
Let $1 \leq p <\infty$ and $\lambda \geq 0.$ 

By  {\bf Campanato spaces} $ {\cal L}^{p,\lambda} (\Omega)$ we denote the linear space of
functions $u \in L^p(\Omega)$ such that 
$$
[f]_{p,\lambda} \,=\, \Biggl\{  \sup_{{x \in \Omega ,}\ {0<\rho < \mathrm{diam}
\,\Omega}} 
\rho^{-\lambda} \int\limits_{\Omega \cap B(x,\rho)} |f(y)- f_\rho |^p
dy \Biggr\}^{\frac{1}{p} }< + \infty.
$$
\end{defin}

$ {\cal L}^{p,\lambda} (\Omega)$ are Banach spaces with the
following norm
$$
\| f \|_{{\cal L}^{p,\lambda} (\Omega)} =  \|f\|_{{L}^{p} (\Omega)} + [f]_{p,\lambda}
$$
which simply demonstrate that $u \in {\cal L}^{p,\lambda} (\Omega)$
if and only if 
$$
\sup_{{x \in \Omega ,} {0<\rho < \mathrm{diam}\,\Omega}} \rho^{-\lambda} 
\inf_{c \in \R} \int_{\Omega \cap B(x,\rho)}|f\,-\,c|^p dy < \infty.
$$

Using H\"older inequality we have that
$$
{\cal L}^{p_1,\lambda_1} (\Omega) \subset {\cal L}^{p,\lambda} (\Omega)
$$
where
$$
p\,\leq \, p_1 ,\qquad \frac{n - \lambda}{p} \geq {\frac{n - \lambda_1}{p_1} }.
$$

Let us observe that
$$
\int_{\Omega \cap B(x,\rho)}|f\,-\,f_\rho|^p dy \leq {\cal C}\cdot 
\int_{\Omega \cap B(x,\rho)}\Bigg( |f|^p  + |\Omega \cap B(x,\rho) | \cdot |f_\rho|^p 
\Bigg) dy 
$$
and also
$$
|f_\rho|^p \leq \frac{1}{|\Omega \cap B(x,\rho)|} \cdot \int_\Omega
|f|^p dy.
$$
Below, we consider $0 \leq \lambda < n.$ 

We use 
$$
[f]_{p,\lambda}  \leq C \|f\|_{L^{p,\lambda} (\Omega)} .
$$
to obtain the following relation between Morrey and Campanato spaces
\begin{equation}
L^{p,\lambda} (\Omega) \subset {\cal L}^{p,\lambda} (\Omega) .
\end{equation}

Let us now recall the definitions of the $BMO$ and $VMO$ classes.
\begin{defin}\label{def2.3}(see \cite{J-N}).
We say that a function $f$ belongs to the John-Nirenberg space $BMO,$ 
or that $f$ has  "bounded mean oscillation", if 
$$
\| f \|_* \equiv
\sup_{{B\rho} \subset {\R}^n} \frac{1}{|B_\rho|}
\int\limits_{B_\rho} |f(y)-f_\rho| dy < \infty
$$
where $f_\rho$ is the integral average 
of the function $f$ over the balls $B_\rho .$ 
\end{defin}

Let us define, for a function $f\in BMO,$
$$
\eta (r) = \sup_{{x \in {\R}^n,}\ {\rho \leq r}}
\frac{1}{|B_\rho|} \int\limits_{B_\rho} |f(x)-f_\rho| dx.
$$

\begin{defin}\label{def2.4}(see \cite{Sarason}).
A function $f\in BMO$ belongs to the class $VMO,$ or $f$ has
"vanishing mean oscillation" if 
$$\lim_{r \to 0^+} \eta (r) =0.$$
\end{defin}

We are now ready to formulate the hypothesis on the terms $A_{ij}^{\alpha \beta}(x,u)$ and $g(\cdot ,u,Du ).$

We suppose that $A^{\alpha \beta}_{i j}(x,u)$ are bounded functions in $\Omega
\times \R^N ,$ such that:
\begin{enumerate}
\item[(A1)]  
     $A^{\alpha\beta}_{ij} = A^{\beta\alpha}_{ji}$.
\item[(A2)]  For every $u \in \R^N$, $A^{\alpha\beta}_{ij}(\cdot , u)
    \in VMO(\Omega)$.
\item[(A3)]  For every $x\in \Omega$ and $u,v\in \R^N$,
    \[ \big|A^{\alpha\beta}_{ij}(x,u)-A^{\alpha\beta}_{ij}(x,v)\big|
       \leq \omega(|u-v|^2)
    \]
    for some monotone increasing concave function $\omega$ with $\omega(0)
    =0$.
\item[(A4)] There exists a positive constant $\nu$
    such that
    \[
       \nu |\xi|^2 \leq A^{\alpha\beta}_{ij}(x,u)\xi^i_\alpha \xi^j_\beta 
    \]
    for a.e. $x\in\Omega$, all $u\in \R^N$ and $\xi \in \R^{nN}$.
\end{enumerate}
We suppose that the function $g$ is a Charath\'eodory function, that is: 
\begin{enumerate}
\item[(g1)]  $g(\cdot ,u,Du )$ is measurable in $x$ $\forall u \in \R^N,
\forall z \in \R^{nN};$
\item[(g2)]  $g(x,\cdot , \cdot )$ is continuous in $(u,z)$ $\, a. \,\,e. \, x \in \Omega;$

moreover we consider $g$ satisfying the condition:
\item[(g3)] 
$$
|g(x,u,z)|\, \leq \, g_1 (x) \,+\, H \,|z|^\gamma ,
$$
$g_1 \geq 0,$ a. e. in $\Omega,$ $ g_1 \in L^p(\Omega), \,$
$2<p \leq \infty,$ $H \geq 0, $ $0 \leq \gamma < 2.$
\end{enumerate}

We point out that since $C^0$ is a proper subset of VMO, the
continuity of $A^{\alpha \beta}_{i j}(x,u)$ with respect to $x$ is not
assumed. 

Let us make some remarks on $VMO$ class.
It was at first defined by Sarason  in 1975 and later it was considered 
 by  many others. For instance we recall 
the papers by Chiarenza, Frasca and Longo \cite{CFL1}  where the authors answer a question raised thirty years
before by C. Miranda in \cite{Mir}. In his paper he considers a linear elliptic equation 
where the coefficients $a_{ij}$ of the higher order derivatives are in
the class $W^{1,n}(\Omega)$  and  asks whether the gradient of the solution is bounded, if $p>n.$ 
Chiarenza, Frasca and Longo suppose $a_{ij} \in VMO$ 
and prove that $Du$ is h\"older continuous for all $p \in ]1, +
\infty[.$
We point out that $W^{1,n} \subset VMO$ because, using Poincare's inequality 
$$
\intmean_B |f(x) -f_B| \leq c(n) \left( \int_B |\nabla u| dx
\right)^{\frac{1}{n}}
$$
and the term on the right-hand side tending to zero as $|B|\to 0.$
Later the interior estimates obtained by Chiarenza, Frasca and Longo were extended 
to boundary estimates in \cite{CFL2}. From these papers on,
many authors have used this space VMO to obtain regularity results
for partial differential equations and systems with discontinuous coefficients. We recall for
example Bramanti and Cerutti \cite{BC} for parabolic
equations  
and many others.

With this useful assumption we investigate the regularity of the minimizers for the
quadratic functional. Its existence is
guaranteed, being the functional ${\cal A}$
sequentially lower semicontinuous with respect to the $H^{1,2}-$weak
topology (see \cite{5}).

\section{Main Results}


\begin{thm} 

Let $u \in W^{1,2}(\Omega, {\R}^N)\,$ be a minimum of the functional
${\cal A}(u, \Omega)$ defined above. 
Suppose that assumptions (A-1), (A-2), (A-3), (A-4), $1<q \leq 2$ and (g-1), (g-2) and (g-3) are satisfied.
Then for $\lambda = n(1- \frac{q}{p})$
we have
\begin{equation}
D\,u \,\in {\cal L}^{q, \lambda}_{\mathrm{loc}}(\Omega_0, {\R}^{n \, N})
\end{equation}
where 
$$
\Omega_0=\{\,  x \in \Omega \colon \liminf_{R \to 0} \!\!\frac
{1}{R^{n-2}}\int_{B(x,R)} \!\!\!|Du(y)|^2 dy =0 \,
\}.
$$
\end{thm}

The set $\Omega_0 $ is obligatory, in fact when we pass from
the regularity theory for scalar minimizers of solutions of elliptic
equations to the regularity theory for vector-valued minimizers of
solutions of elliptic systems, the situation changes completely: 
regularity is an exceptional occurence everywhere, excluding the two dimensional case. In 1968 De Giorgi in \cite{2} showed that his regularity result for solutions of second
order elliptic equations with measurable bounded coefficients cannot
be extended to solutions of elliptic systems.
He presented the quadratic functional
$$
{\cal S}\,=\, 
\int_\Omega A^{\alpha \beta}_{i j}(x) D_{\alpha} u^i D_{\beta} u^j d x 
$$ 
with $A^{\alpha \beta}_{i j} \in L^\infty(\Omega), $ such that
$$
\!\!\!\!\!\exists  \nu \! > \! 0  \colon \! A^{\alpha \beta}_{i
j}\! \chi_{\alpha}^{i}\! \chi_{\beta}^{j} \geq \nu |\chi |^2 ,\!\! \!\!\!\!
 ~~\mathrm{a.e.}\!\!~ x \in \Omega,\!\! \!\!~\forall \!\chi \in  \!{\R}^{n N} \!\!\!\!.
$$ 
De Giorgi proves that ${\cal S}$ has a minimizer that is a function having
a point of discontinuity in the origin.   
Later, Sou\v{c}ek in \cite{Soucek} showed that minimizers of functionals of the type
${\cal S}$ can be discontinuous, not only in a point, but also on a dense subset of $\Omega.$
Modifying De Giorgi's example, Giusti and Miranda in
\cite{Giusti-Miranda2}  showed that solutions of elliptic quasilinear systems of the type
$$
\int_\Omega A^{\alpha \beta}_{i j}(u) D_{\alpha} u^i D_{\beta} 
\varphi^j d x\,=\, 0, \forall \varphi \in C^\infty_0 (\Omega, {\R}^{N} )
$$
with analytic elliptic coefficients $A^{\alpha \beta}_{i j}$ have
singularities in dimension $n \geq 3.$ 
We observe that we can get global regularity for some special cases,
see for example \cite{Tachi1}.

Similar examples were presented in the meantime independently by
Maz'ya in \cite{Maz}.
Even Giaquinta in 1993 in \cite{5}, Morrey in \cite{Morrey-book} and others were
interested in these problems of solutions of elliptic systems,
solutions in general non regular.
Then we can prove regularity  except on a set, hopefully not too
large. 

For linear systems, regularity results assuming $A^{\alpha \beta}_{i
j}$ constant or in $C^0(\Omega),$ have been obtained by  Campanato $\,\,$
 (see \cite{1}).
Without assuming continuity of coefficients, we mention the study by
Acquistapace \cite{Acq} where Campanato's results are refined
considering that $A^{\alpha \beta}_{i j}$ belongs to a class that neither contains nor is being contained by $C^0(\Omega),$ 
hence in general discontinuous.
Moreover, we recall the study made by Huang in \cite{QH} where he shows the regularity of weak solutions of linear
elliptic systems with coefficients $A^{\alpha \beta}_{i j}(x) \in
VMO.$ So, it seems to be natural to expect partial regularity
results under the condition that the coefficients of the principal
terms $A^{\alpha \beta}_{i j} \in VMO,$  even for nonlinear cases.
Dane\v{c}\v{e}k and Viszus in \cite{3} 
consider the
regularity of minimizer for the functional
$$
\int_\Omega  \bigl\{
A^{\alpha \beta}_{i j}(x) D_{\alpha} u^i D_{\beta} u^j + g(x,u,Du)  
\bigr\} dx
$$
where the term $g(x,u,Du)$ is such that
$$
|g(x,u,z)|\, \leq \, f(x) \,+\, |z|^\gamma
$$
where $f \in L^p(\Omega), $ $2<p\leq \infty , $ $f \geq 0 $ a. e. on
$\Omega,$  $L$ is a non-negative constant and $0 \leq \gamma <2.$

They obtained h\"older regularity of minimizer assuming that $A^{\alpha
\beta}_{i j}(x) \in VMO.$

We also recall the paper \cite{DES} where the authors obtain regularity results for minimizers
of the quasilinear functionals 
\[
      \int_\Omega 
      A^{\alpha\beta}_{ij}(x, u) D_\alpha u^i D_\beta u^j dx,
\]
where the coefficients $A^{\alpha\beta}_{ij}(x,u)$ have  VMO dependence on
the variable $x$ and continuous dependence on $u.$

In the paper \cite{RT1} we improve the last mentioned result in Morrey spaces 
because we have considered inside the integral the term $g(x,u,Du)$ and the result by Dane\v{c}\v{e}k and 
Viszus because we consider $A^{\alpha \beta}_{i j}$ dependent not only on $x$ but also on $u.$ 
In the present note we extend our above cited
regularity results because we consider the more general class Campanato spaces.

Before an outline of the proof we state a preliminary Lemma by
Campanato. 

\begin{lem} 
Let $B(x_0,R)$be a fixed ball and $u \!\in \!W^{1,2}(B(x_0,R)\!; \!\R^N)$ be a weak
solution of the system
$$
D_\alpha (A^{\alpha \beta}_{i j}D_{\beta} u^j ) = 0, \quad i=1,\ldots,N
$$
where $A^{\alpha \beta}_{i j}$ are constant and satisfy the
ellipticity condition.  Then $\forall t \in (0,1]$ 
$$
\int_{B(x_0, tR)}|Du|^2 dx \leq c \cdot t^n \int_{B(x_0,R)}|Du|^2 dx .
$$
\end{lem}

{\it PROOF OF THE THEOREM }

For simplicity we 'll consider the case $g=0.$ 
Let $R>0$ and $x_0 \in \Omega$ such that $B(x_0,R)\subset
\subset \Omega.$

Let $v$ be the minimum of the ``freezing'' functional ${\cal A},$ that
is 
$$
{\cal A}_0(v, B(x_0, \frac{R}{2})) \,=\, \int_{ B(x_0, \frac{R}{2}) }
A^{\alpha \beta}_{i j}(x_0,u_{\frac{R}{2}}) D_{\alpha} v^i D_{\beta} v^j dx
$$
with $v\equiv u$ on $\partial B(x_0 \frac{R}{2}).$ 
The idea of freezing is the same used by Chiarenza, Frasca and Longo in \cite{CFL1}.

For $0\leq \lambda <n$ and $q \leq 2$ we have that

$$ 
\|Du\|_{ {\cal L}^{q,\lambda} (\Omega)}= \|Du\|_{{ L}^{q}
(\Omega)}+ [Du]_{q,\lambda}\leq 
$$ 
$$
\,
$$ 
$$
\qquad\qquad\qquad\qquad\qquad \leq {\cal K} \,\|Du\|_{{
L}^{q,\lambda} (\Omega)}\,\leq {\cal K} \,\|Du\|_{{ L}^{2,\lambda}
(\Omega)}
$$ 
where the constant ${\cal K}$ is independent of $u.$
We observe that $A^{\alpha \beta}_{i j}(x_0,u_{\frac{R}{2}})$ are constant
coefficients, then from the above Lemma, $\forall t \in (0,1],$ $$
\!\!\!\!\!\int_{B(x_0, t\frac{R}{2} )}|Du|^2 dx \leq \!c \!\cdot t^n
\int_{B(x_0,\frac{R}{2} )}|Du|^2 dx .  $$ Let $w= u-v,$ then $w \in
W^{1,2}_0( B(x_0,\frac{R}{2}))$

$$
\begin{array}{cl}
   & \displaystyle{\int_{B(x_0,\frac{t\,R}{2})} |Du|^2 dx }\leq \\
   & \\
   \leq & \displaystyle{c \cdot 
   \left\{  t^n 
   \int_{B(x_0,\frac{R}{2})} |Du|^2 dx +
   \int_{B(x_0,\frac{R}{2})} |Dw |^2 dx
\right\}
.}
\end{array}
$$
Let us estimate:
$$
\int_{B(x_0,\frac{R}{2})} \,\,|Dw |^2 \,dx.\,
$$

From the hypothesis and a Lemma in \cite{9} we have
$$
\!\!\!\!\nu \!\!\int_{B(x_0,\frac{R}{2})} \!\!\!|Dw |^2 dx \!\!
\leq  \!\!
\left\{ \!
{\cal A}^0 (u,B(x_0, \!\frac{R}{2})) \!- \! \!\!
{\cal A}^0 (v,B(x_0,\! \frac{R}{2})) 
\right\}
$$
adding and subtracting:
$$
A^{\alpha \beta}_{i j}(x,u_{\frac{R}{2}}) D_{\alpha} u^i D_{\beta} u^j, \quad 
A^{\alpha \beta}_{i j}(x,u_{\frac{R}{2}}) D_{\alpha} v^i D_{\beta} v^j 
$$
$$
A^{\alpha \beta}_{i j}(x,u) D_{\alpha} u^i D_{\beta} u^j\,, \quad 
A^{\alpha \beta}_{i j}(x,v) \,D_{\alpha} v^i D_{\beta} v^j 
$$
we obtain 
some different kinds  of integrals that we now examine.

Using $L^p$ estimate, we can estimate the terms with
$|Du|^2$ as follows: 
$$
\int_{B(x_0,\frac{R}{2})}  |A^{\alpha \beta}_{i
j}(x_0,u_{\frac{R}{2}})-A^{\alpha \beta}_{i j}(x,u_{\frac{R}{2}})|\, \cdot \,|Du|^2\,dx\,\leq 
$$
$$
\leq c\,\big\{\eta (A(\cdot,u_{\frac{R}{2}});R)\big\}^{1-\frac{2}{p} }  \,\int_{B(x_0,R)}  |\,Du\, |^2\,dx.
$$

We also observe that: 
$$
\int_{B(x_0,\frac{R}{2})}  | A^{\alpha \beta}_{i j}(x,u_{\frac{R}{2}})-A^{\alpha
\beta}_{i j}(x,u)|\,D_{\alpha} u^i D_{\beta} u^j \, \,dx\,\leq 
$$
$$
\leq c
\left(
\int_{B(x_0,R)} \!\! |Du|^2 dx \!
\right)\!\! 
\left( \!\!
\intmean_{B(x_0,\frac{R}{2})}  \!\!\omega(|u_{\frac{R}{2}}-u|^2) dx \!\!
\right)^{ \!\!1\!-\!\frac{2}{p} } \!\leq  
$$
$$
\leq c
\left(\int_{B(x_0,R)} \!\! |Du|^2 dx \!\right)
    \left(\omega 
    \left( R^{2-n}
\int_{B(x_0,\frac{R}{2})} |D\,u|^2\,dx\, \right) \right)^{ 1-\frac{2}{p}}.
$$

Moreover, we can estimate the terms having $|Dv|^2$ similarly
using $L^p$ estimates for $Dv$.

Then if $\rho = t R$
$$
\int_{B(x_0,\rho )} |Du|^2dx \leq 
$$
$$
\leq  c \Biggl\{ \left( \frac{\rho}{R} \right)^n+
\left(\omega
\left(R^{2-n}
\int_{B(x_0,\frac{R}{2})} |D\,u|^2\,dx\,
\right)
\right)^{ 1-\frac{2}{p}}
$$
$$
+\biggl( \eta (A(\cdot,u_{\frac{R}{2}}),R) \biggr)^{ 1-\frac{2}{p}} \Biggr\}
 \cdot\left(\int_{B(x_0,\frac{R}{2})} 
|D\,u|^2 dx\right).
$$
Using a Lemma contained in \cite{4} and selecting $\rho$ sufficiently
small, specifically $\rho <\frac{R}{2},$ we have
$$
\int_{B(x_0,\rho )} |Du|^2dx \leq c \cdot \rho^\lambda .
$$


M.A. Ragusa,
	 Dipartimento di Matematica, Universit\`{a} di Catania,
	 Viale A. Doria, 6, 95125 Catania, Italia
	 
\it E-mail address: \texttt{maragusa@dipmat.unict.it}\\
\rm
~~\\
A. Tachikawa,
     Department of Mathematics, Faculty of Science and Technology,
     Tokyo University of Science,
     Noda, Chiba, 278-8510, Japan
     
\it E-mail address: \texttt{tachikawa$ \_ $atsushi@ma.noda.tus.ac.jp}
\end{document}